\documentclass[reqno,11pt]{amsart}
\usepackage{amssymb,latexsym}
\usepackage{hyperref}
\usepackage{color}

\newtheorem{thm}{Theorem}
\newtheorem{lemma}{Lemma}

\newtheorem{prop}{Proposition}

\theoremstyle{definition}

\theoremstyle{remark}

\makeatletter
\@addtoreset{figure}{section}
\def\thefigure{\thesection.\@arabic\c@figure}
\def\fps@figure{h, t}
\@addtoreset{table}{bsection}
\def\thetable{\thesection.\@arabic\c@table}
\def\fps@table{h, t}
\@addtoreset{equation}{section}

\makeatother

\allowdisplaybreaks
\def\intprod{\mathbin{\hbox to 6pt{%
                 \vrule height0.4pt width5pt depth0pt
                 \kern-.4pt
                 \vrule height6pt width0.4pt depth0pt\hss}}}

\begin{document}

\title[Liquid crystals on Riemannian manifolds]
{Well-posedness and Global Attractors for Liquid crystals on 
Riemannian manifolds}

\author{Steve Shkoller}

\address{Department of Mathematics\\
University of California \\
Davis, CA 95616}
\email{shkoller@math.ucdavis.edu}


\date{Oct. 1, 2000; current version \today} 

\begin{abstract}
We study the coupled Navier-Stokes Ginzburg-Landau model of nematic
liquid crystals introduced by F.H. Lin, which is a simplified version of
the Ericksen-Leslie system.   We generalize the model to
compact $n$-dimensional Riemannian manifolds, and show that the system comes 
from a variational principle.  We present a new simple proof for the local
well-posedness of this coupled system without using the higher-order 
energy law.  We then prove that this system is globally well-posed 
and has compact global attractors when the dimension of the manifold $M$ is two.
Finally, we introduce the Lagrangian averaged liquid crystal equations,
which arise from averaging the Navier-Stokes fluid motion over small spatial
scales in the variational principle.  We show that this averaged system 
is globally well-posed and has compact global attractors even when
$M$ is three-dimensional.

\end{abstract}

\maketitle

\tableofcontents

\section{Introduction}
\label{Intro}

Nematic liquid crystals are well-studied and interesting examples of 
anisotropic non-Newtonian fluids.  A liquid crystal is a phase of a material
between the solid and liquid phases.  The solid phase has strong 
intermolecular forces that keep the molecular position and orientation fixed, 
while in the
liquid phase, the molecules neither occupy a specific average position 
nor do they remain in any particular orientation; the liquid crystal phase
does not have any positional order, but does possess a certain amount of
orientational order.
This phase is described by a velocity field,
as well as a director field that describes locally the averaged direction or
orientation
of the constituent molecules.  In this paper, we shall analyze the behavior
of a certain model of nematic liquid crystals on compact Riemannian manifolds.

We let $(M,g)$ denote a 
smooth, compact, connected, $n$-dimensional Riemannian manifold with 
smooth (possibly empty) boundary ${\partial M}$.  If ${\partial M}=
\emptyset$, then we assume that the Euler characteristic $\chi(M)$ does not
vanish.  We study the following system of nonlinear partial differential
equations:
\begin{subequations}
  \label{lc}
\begin{gather}
u_t + \nabla_uu = -\operatorname{grad}p + \nu \operatorname{Div}
\operatorname{Def}u - \lambda\operatorname{Div}(\nabla d^T\cdot \nabla d) \,,
         \label{lc.a}\\
 \operatorname{div}u(t,x)=0 \,,
         \label{lc.b}\\
d_t + \nabla_ud = \gamma\left(\hat\Delta d - 
\frac{1}{\epsilon^2}(|d|^2 -1)d\right) \,,
         \label{lc.c}\\
u=0 \ \rm{on} \ {\partial M}, \ \ d = h \ \rm{on} \ {\partial M} \ \ \
g(h,h)=1 \, \rm{\text{ or }} {\partial M} = \emptyset\,,
         \label{lc.d}\\
u(0,x)=u_0, \ \ d(0,x) = d_0 \text{ and } d_0|_{\partial M}=h \text{ if }
   {\partial M} \neq \emptyset .
         \label{lc.e}
\end{gather}
\end{subequations}
Here $u(t,x)$ and $d(t,x)$ are time-dependent vector fields on $M$,
$\nabla$ denotes the Levi-Civita covariant derivative associated to the
Riemannian metric $g$, 
$\operatorname{Def} u = {\frac{1}{2}}(\nabla u + \nabla u^T)$ denotes 
the (rate of) deformation tensor, the superscript $(\cdot)^T$ denotes the
transpose, $\hat\Delta$ denotes the rough Laplacian of $g$ defined in 
(\ref{rough}),
and $\nu, \lambda, \gamma$ are positive constants.
In the case that $M$ is flat, an open subset of Euclidean space for instance,
then $\nabla$
is the componentwise gradient,  and $\hat\Delta$ is the componentwise
Laplacian given in coordinates $x^i$ by $\hat\Delta = \sum_{i=1}^n \frac
{\partial^2}{\partial x^i \partial x^i}$.  The system
of equations (\ref{lc}) is the simplified Ericksen-Leslie model 
\cite{E1,E2,L} of nematic liquid
crystals first introduced by F.H. Lin in \cite{Lin1} and later analyzed
by F.H. Lin and C. Liu in \cite{LL1,LL2}.   Beautiful numerical simulations
can be found in \cite{LW}.

This system couples the Navier-Stokes (NS) equations with the Ginzburg-Landau
(GL)
penalization of the harmonic map heat flow.  The vector field $u(t,x)$ is
the velocity field of the fluid, while $d(t,x)$ is the penalized 
(Ginzburg-Landau) approximation to the unit-length director field, representing
the orientation parameter of the nematic liquid-crystal.   

 The parameter
$\epsilon>0$ is the penalization parameter, $\nu$ denotes the kinematic
viscosity of the fluid, $\lambda$ is an elastic constant, and $\gamma$
is the relaxation-time parameter.

The coupling term $\operatorname{Div}(\nabla d^T \cdot \nabla d)$ 
preserves the regularity of the Navier-Stokes equations:  when velocity
$u=0$, (\ref{lc.a}) becomes $\operatorname{Div}(\nabla d^T \cdot \nabla d)
= -\operatorname{grad} p$,  so that when $d_t=0$, $d$ is a solution
of the static portion of (\ref{lc.c})  together with the constraint
$\operatorname{Div}(\nabla d^T \cdot \nabla d)
= -\operatorname{grad} p$.  This constraint pushes the gradient flow
towards a ``regular'' ($H^s$, $s$ sufficiently large)  stationary solution.
Even though the coupling term has two derivatives, analytically, it is
essentially identical to the advection term $\nabla _uu$.

\medskip \noindent
{\bf Results.}  We begin by extending the simplified Ericksen-Leslie model 
of F.H. Lin to a 
compact Riemannian manifold with boundary.  This extension introduces a new 
curvature term in the basic energy laws, and provides a covariant 
(coordinate-independent) description of the liquid crystal dynamics.  
Motion on the sphere is an important application.

In Section \ref{sec_variational}, we prove that the system of equations
(\ref{lc}) actually arises from a simple variational principle 
(the system was originally derived using balance laws).  The variational
principle is the key to our analysis, for it gives the correct scaling;
namely, it shows that when $d$ is taken to have  one derivative greater
regularity than $u$, the liquid crystal system behaves as if it were parabolic,
a fact which was not previously known (see \cite{LL1}).  In fact, according
to \cite{Liu}, because the interaction term 
$\operatorname{Div}(\nabla d^T \cdot \nabla d)$ 
formally has as many derivatives as the diffusion term, the standard 
Galerkin procedure for obtaining local solutions had failed in prior attempts.

In Section \ref{sec_local}, we give a very simple proof of local well-posedness
of the system (\ref{lc})  on $(M,g)$ (in Theorems \ref{thm1} and \ref{thm1a})
using the contraction mapping theorem;  
this significantly simplifies the clever, but lengthy, {\it modified Galerkin}
procedure employed by Lin and Liu in \cite{LL1}.  Moreover, the proof does
not require use of either the maximum principle or higher-order energy laws.

In Section \ref{sec_energy}, 
we show (in Propositions \ref{prop1}, \ref{prop1a}, and \ref{H1H2})
that on two-dimensional Riemannian manifolds with smooth boundary
(possibly empty), there exists  an absorbing set for $u$ in $H^1$ and $d$ in 
$H^2$.  

In Section \ref{sec_global}, we prove the global well-posedness of
the system (\ref{lc}), as well as the existence of absorbing sets for $u$ in 
$H^s$ and $d$ in $H^{s+1}$, and hence of a compact global attractor when the 
dimension is $n=2$ (see Theorems \ref{thm2} and \ref{thm2a}).  
When ${\partial M}=\emptyset$ brute-force energy estimates
may be computed, but when ${\partial M}$ is not empty, we use the
Ladyzhenskaya method to obtain the uniform bounds.  We remark that the existence
of global attractors for this system was not previously known.  We also 
remark that since the Navier-Stokes equations are a subsystem of 
(\ref{lc}), one does not expect to be able to prove results in dimension
three which do not already exist for the Navier-Stokes equations; namely,
the problem of unique classical solutions remains open, while weak
solutions exist \cite{LL1}.

Finally, in Section \ref{sec_averaged}, we introduce the Lagrangian averaged
liquid crystal equations (\ref{avg_lc}).  
This system is based on the Lagrangian averaged
Navier-Stokes equations (see \cite{MS1} and references therein), and is
derived by averaging the Navier-Stokes flow over small spatial scales which
are smaller than some positive small number $\alpha$.  We show that this
averaged system retains the structure of the original system derived by
Lin in the form of averaged energy laws, but has the advantage of being
globally well-posed on three-dimensional domains (see Theorem \ref{thm3}).  
The averaged energy
law shows that when both the fluid flow is averaged together with the director
field, both $u$ and $d$ scale similarly, and $d$ is not required to have
one-derivative greater regularity.  Of course, physically, it seems much
more natural to us to average the fluid flow, since the molecular orientation
is already an averaged quantity.  We believe that the averaged liquid crystal
system will be the ideal model for numerical computation.

\bigskip
\noindent
{\bf Some Notation and Interpolation Inequalities.}
We shall use the notation $H^s(TM)$ to denote the $H^s$-class vector
fields on the manifold $M$.
The $H^s(TM)$ inner-product is given, in any local chart, by
$$\langle u, v \rangle_s = \sum_{|\alpha|=0}^s \langle D^\alpha u,
D^\alpha v\rangle,$$
where 
$$\langle u, v \rangle = \int_M g(x)\left( u(x), v(x) \right) \mu(x)$$
denotes the $L^2$ inner-product,  $\alpha=(\alpha_1, ..., \alpha_n)$
is a multi-index, and 
$$|\alpha| = \alpha_1 + \cdot\cdot\cdot +\alpha_n, \ \ \ \
D^\alpha = \partial_1^{\alpha_1} \cdot\cdot\cdot \partial_n^{\alpha_n}.$$
We shall denote the $H^s(TM)$ norm 
by 
$$|u|_s = \langle u, u \rangle_s,$$
$\langle \cdot, \cdot \rangle = \langle \cdot, \cdot \rangle_0$,
and $| \cdot | = |\cdot |_0$.  We set $H^1_0(TM)$ to consist of
those vector fields in $H^1(TM)$ which have zero trace on ${\partial M}$.
Similarly, vectors in $H^1_h(TM)$ have trace $h$ on ${\partial M}$.
We let
$$H^s_h(TM) = H^s(TM)\cap H^{s-{\frac{1}{2}}}(T{\partial M}), \ \
g(x)(h(x),h(x))=1 \ \forall x\in{\partial M}\ \ s\ge 1 $$
denote the space of $H^s$ vector fields on $M$ which have ($H^1$) trace $h$
on ${\partial M}$ and where $h\in H^{s-{\frac{1}{2}}}(T{\partial M})$.

For each $x\in M$, we let 
$B^\delta_x =\{v\in T_xM \ | \ g(x)(v,v) \le \delta\}$,
and set $B^\delta = \cup_{x\in M} B^\delta_x$.  We let
$H^s(M,B^\delta)$ denote the $H^s$-class maps from $M$ into $B^\delta$.

We have the product rule
$$ D^\alpha( f \ g) = \sum_{ \stackrel{|\beta|\le|\alpha|}{\alpha -\beta>0} }
c_{\alpha,\beta} \left(D^\beta f\right) \ \left(D^{\alpha - \beta} g\right).$$
For any integer $s\ge 0$, we set
$$D^su=\{D^\alpha u \ : \ |\alpha|=s\}, \ \ \ \ \ \ \ 
\|D^su\|_{L^p}= \sum_{|\alpha|=s} \| D^\alpha u\|_{L^p}.$$

We define the spaces
\begin{align*}
{\mathcal V}& =\{ u \in C^\infty(TM) \ | \ \operatorname{div}u=0, \
\ g(u,n)=0 \ \rm{on} \ {\partial M}\}, \\
{\mathcal W}& =\{ u \in C^\infty_0(TM) \ | \ \operatorname{div}u=0\},
\end{align*}
and through-out the paper, we shall use $W^s$ and $V^s$ denote the
closure in $H^s$ of ${\mathcal V}$ and ${\mathcal W}$, respectively.  It follows
that
\begin{align*}
V^s &= \{ u\in H^s(TM) \ | \ \operatorname{div}u=0,\ 
g(u,n)|_{{\partial M}}=0\},\\
W^s &= \{ u\in H^s(TM) \cap H^1_0(TM) \ | \ \operatorname{div}u=0\}.
\end{align*}
In section \ref{sec_local}, we shall give an equivalent definition of
$W^s$ using powers of the Stokes operator.

We shall need some standard interpolation inequalities, which follow from
the {\it Gagliardo-Nirenberg inequalities } \cite{N},\cite{T}:

{\it Suppose}
\[
\frac{1}{p} = \frac{i}{n} + a \left( \frac{1}{r} - \frac{m}{n}
\right) + (1 - a ) \frac{1}{q}
\]
{\it where $ i/m \leq a \leq 1$ \textup{(}if $m - i - n/r $ is
an integer $ \geq 1 $, only $ a < 1 $ is allowed\textup{)}.
Then for $ f : M\rightarrow TM $},
\begin{equation} \label{gagnir}
|D^if|_{L^{p}} \leq C|D^m f|^a_{L^{r}}\cdot
|f|^{1-a}_{L^{q}}
\end{equation}
In what follows, we
shall use $C$ as a generic constant.  Some specific cases
in two dimensions ($n = 2 $) that we shall need
are as follows:

\begin{align}
| v | _{L ^{\infty} } & \leq C | D ^2 v | _{L^2} ^{1/2}
|v|_{L^2} ^{1/2} \label{interp1}\\
|v| _{L ^4} & \leq C | D v  | _{L^2}^{1/2} |v|
_{L^2}^{1/2}
\label{interp2}\\ | D ^i v  | _{ L ^2} & \leq C |v| _{L ^2} ^{1
- i / m } | D ^m v | _{ L ^2} ^{ i/m }. \label{interp3}
\end{align}

Equation (\ref{interp1}) is often called the Agmon inequality, while
(\ref{interp2})-(\ref{interp3}) are often referred to as the 
Ladyzhenskaya inequalities.

We shall use $\operatorname{div}$ for the divergence operator on vector
fields, and $\operatorname{Div}$ for the divergence operator on sections of
$T^*M \otimes TM$.

\section{The Variational Principle}
\label{sec_variational}
In this section, we shall explain how the system of equations (\ref{lc}) 
arise from a simple variational principle, (\ref{lc.a}) being the first 
variation of the action with respect to the Lagrangian flow variable, and
(\ref{lc.c}) being the $L^2$ gradient flow of the first variation of
the action with respect to the director field.  It was not previously
known that (\ref{lc}) can be obtained from a variational principle; rather,
balance arguments were invoked to derive the model.

We let $\eta(t,x)$ denote the Lagrangian flow variable, a solution of the
differential equation
$$ \partial_t \eta(t,x) = u(t,\eta(t,x), \ \ \ \eta(0,x) =x.$$
For $I=[0,T]$, and each $t\in I$, $u\in C^0(I,W^s)$, $s>(n/2)+2$, the map 
$\eta(t,\cdot): M \rightarrow M$ is an $H^s$ volume-preserving
diffeomorphism with $H^s$ inverse, and restricts to the identity map 
on the boundary ${\partial M}$.  We shall denote
this set of maps by ${\mathcal D}_{\mu,D}^s$.  It is a fact, that
for $s>(n/2)+1$, the set ${\mathcal D}_{\mu,D}^s$ is a $C^\infty$ 
(weak) Riemannian manifold (see \cite{EM} and \cite{S}).

We define the action function  $S: {\mathcal D}_{\mu,D}^s \times H^{s+1}(TM)
\cap H^1_0(TM) \rightarrow {\mathbb R}$ by
\begin{align}
&S(\eta, d) =  {\frac{1}{2}}\int_I \int_M \left\{
g(\eta(x)) \left( u(t,\eta(t,x)), u(t,\eta(t,x))\right)\right. \nonumber\\
&\qquad\qquad \left. + \lambda
g(\eta(x)) \left(
\nabla\left[ d(t,\eta(t,x)) \right],
\nabla\left[ d(t,\eta(t,x)) \right] \right) + 2 F (d)
\right\} \mu\, dt,
\label{action}
\end{align}
where $F(d) = \frac{1}{4\epsilon^2}(|d|^2 -1)^2$.  Notice that
$$ f(d) = \operatorname{grad} F(d),$$
where 
$$f(d) \equiv \frac{1}{\epsilon^2}\left(|d|^2 -1\right)d$$
is the (GL) nonlinearity  in (\ref{lc.c}).
The first term on the right-hand-side of (\ref{action})
is the {\it kinetic energy} of the fluid, the second term is the
{\it elastic energy} of the polymers, and the third term is the
{\it unit-length constraint} on the director field $d$.  
As a consequence of the right-invariance of $S$ with respect to the lifted
action of
${\mathcal D}_{\mu,D}^s$, we may compute the kinetic
energy of the fluid as well as the elastic energy 
along the particle trajectory $\eta(t,x)$.  The interaction, or coupling,
between the velocity $u$ and the director $d$ comes precisely from the
elastic energy being computed along the Lagrangian flow $\eta(t,x)$.

The elastic energy $(1/2)\int_M |\nabla d|^2 \mu$ is a simplified form
of the Oseen-Frank energy, given upto the null-Lagrangian by
\begin{equation}\label{frank}
\int_M \left[ \kappa_1 |\operatorname{div} d|^2 + \kappa_2|d \times
\operatorname{curl}d|^2 + \kappa_3 |d \cdot \operatorname{curl} d|^2\right] \mu.
\end{equation}
The terms in the integrand represent, respectively, the energy due to
{\it splay}, {\it bending}, and {\it twisting} of the polymers in the
nematic liquid crystal.  When $\kappa= \kappa_1=\kappa_2=\kappa_3$, then
(\ref{frank}) reduces to $\kappa\int_M |\nabla d|^2 \mu$.
We see that in the Eulerian frame, for  a director field which is exactly
taking values in the unit sphere, the energy is given by
$$
\text{Energy} = {\frac{1}{2}} \int_M\left( |u(x)|^2+ \lambda |\nabla d|^2\right)
\mu.
$$
The penalized form of this energy is then
\begin{equation}
\label{energy}
E = {\frac{1}{2}} \int_M\left( |u(x)|^2+ \lambda |\nabla d|^2
+2 \lambda F(d) \right)
\mu,
\end{equation}
where we suppress the explicit dependence on the small parameter $\epsilon>0$.

The action (\ref{action}) is the right-translated time-integral of the
energy  function
(\ref{energy}).  The penalization was motivated
by the study of harmonic maps of simply-connected domains $\Omega$ into 
spheres (see \cite{BBH}); in particular, the space 
$H^1_h(\Omega,{\mathbb S}^1) = \emptyset$
 when $|\operatorname{degree}(h)| \ge 1$ so that only infinite energy
minimizers exist. As a fix for this problem, the penalization method is
invoked, which
enlarges the space of potential minimizers to $H^1_h(\Omega,{\mathbb R}^2)$ 
(which  is obviously not empty) and simultaneously  imposes the unit-length 
constraint.  

\medskip

To compute the first variation of $S$ with respect to $\eta$, we let
$\varepsilon \mapsto \phi_\varepsilon$ be a smooth curve in 
${\mathcal D}_{\mu,D}^s$ such that $\phi_0 = e$, and
$(d/d \varepsilon)|_{\varepsilon=0}\phi_\varepsilon = w \in W^s$.  
We let $D/d \varepsilon$ denote the covariant derivative along the curve 
$\phi_\varepsilon$.  
Let $\eta^\varepsilon = \eta \circ \phi^\varepsilon$ so that
$$
\eta^0 = \eta, \text{  and  } (d/d \varepsilon)|_{\varepsilon=0} 
\eta^\varepsilon = w \circ \eta.
$$
Then (setting $\lambda=1$ for the moment),
\begin{align*}
&\langle D_1S(\eta,d) , \delta\eta \rangle \\
& \qquad  =
\left.\frac{d}{d \varepsilon}\right|_{\varepsilon=0} S(\eta^\varepsilon,d)
 =
\int_I\int_M  \left\{
g(\eta(t,x))\left( (D/ d \varepsilon)_{\varepsilon =0}
\partial_t \eta^\varepsilon(t,x), \partial_t \eta(t,x) \right) \right.\\
& \qquad
+ \left. g(\eta(t,x))\left( (D/ d \varepsilon)_{\varepsilon =0}
\nabla  \left( d(t,\phi^\varepsilon(\eta(t,x)) \right),
\nabla  \left( d(t,\eta(t,x) \right)
\right) \right\} dx\, dt,
\end{align*}
where $\nabla$ is computed with respect to the moving Lagrangian coordinate
$y=\eta(t,x)$, and where we have used $dx$ to denote the Riemannian
volume-form $\mu$.  We use $D_1$ and $D_2$ to denote the Frech\'{e}t
derivatives of $S$ with respect to $\eta$ and $d$, respectively.  
Integrating by parts, and using the fact that
$\partial_t \eta = u \circ \eta$ and that $\eta$ has Jacobian determinant
equal to one, we see that
\begin{align*}
& \langle D_1S(\eta,d), \delta \eta \rangle
 = \int_I\int_M g(\eta(t,x)) \left(
-((D/dt) \partial_t \eta(t,x), w(t,\eta(t,x) \right) dx \, dt \\
&\qquad + 
 \int_I\int_M g(y) \left( (D/d \varepsilon)|_{\varepsilon=0}
\nabla d \circ \phi_\varepsilon(y) \cdot T\phi_\varepsilon(y),
\nabla d(t,y)\right) dy\, dt \\
&\qquad =
 \int_I\int_M g(y) \left( (-u_t(t,y) - \nabla_uu(t,y) -\operatorname{grad}
p(t,y), w(t,y) \right) dy\, dt \\
&\qquad +
 \int_I\int_M \left\{g(y) \left(  \nabla_w(\nabla d),
\nabla d\right) + g(y)\left( \nabla d(t,y)\cdot \nabla w, \nabla d
\right) \right\} dy\, dt \\
&\qquad =
 \int_I\int_M g(y) \left( (-u_t(t,y) - \nabla_uu(t,y) -\operatorname{grad}
p(t,y), w(t,y) \right) dy\, dt \\
&\qquad +
 \int_I\int_M g(y) \left(  -\operatorname{Div}(\nabla d^T \cdot \nabla d),
w \right) dy\, dt,
\end{align*}
where the last equality follows from the fact that
$\langle \nabla_w(\nabla d), \nabla d\rangle =0$, since 
$\operatorname{div} w =0$.
Thus, since $w$ is an arbitrary variation of $\eta$, we arrive at the
Euler-Lagrange equation
$$u_t + \nabla_uu = -\operatorname{grad} p - \operatorname{Div} 
(\nabla d ^T \cdot \nabla d).$$  The viscosity (diffusion) term follows from
the Ito formula by allowing $\eta(t,x)$ to be a stochastic process, and
replacing deterministic time derivatives with stochastic backward-in-time
mean derivatives (see \cite{G}).  Thus (\ref{lc.a}) follows as the
first variation of the action function $S$ with respect to $\eta$.
Equation (\ref{lc.b}) follows immediately from the fact that $\eta$ is
volume-preserving.

Letting $d^\varepsilon = d+ \varepsilon \delta d$, 
a much simpler computation verifies that 
$$ \langle D_2S(\eta,d) , \delta d \rangle =
\left.\frac{d}{d \varepsilon}\right|_{\varepsilon=0} S(\eta,d^\varepsilon)
 = \int_I\int_M g(y)\left( \hat\Delta d - f(d), \delta d\right) dy\, dt,
$$
where 
\begin{equation}\label{rough}
\hat \Delta d = \nabla^* \nabla
\end{equation}
is the rough Laplacian and $\nabla^*$
is the $L^2$ formal adjoint of the covariant derivative $\nabla$.
Hence, equation (\ref{lc.c}) is simply the $L^2$ gradient flow of
$d\mapsto S(\eta, d)$ given by
$$\frac{d}{dt}(d(t,\eta(t,x)) = D_2S(\eta,d)=\hat\Delta d-f(d).$$

We remark that
\begin{equation}\label{R}
\operatorname{Div}(\nabla d^T \cdot \nabla d) = \nabla d^T \cdot\hat\Delta d
+ g(R(e_i, \cdot)d, \nabla_{e_i}d),
\end{equation}
where $R$ is the Riemannian curvature tensor which is defined for
vector fields $X,Y,Z$ on $M$ by
$$R(X,Y)Z = \nabla_X \nabla_Y Z - \nabla_Y\nabla_XZ + \nabla_{[X,Y]}Z,$$
and where $\{e_i\}$ is any local orthonormal frame.  The curvature term
in equation (\ref{R}) will play an important role in the energy behavior of
the system. 

\section{Local Well-posedness}
\label{sec_local}

Let $P$ denote the Leray orthogonal projection from $L^2(TM)$ onto $W^0$,
and let 
$$A=-P\operatorname{Div}\operatorname{Def}$$ 
denote the Stokes operator, an unbounded, positive,
self-adjoint operator on $W^0$, with domain $D(A)= H^2(TM) \cap W^1$.  As usual,
we set 
$$W^s = D(A^{\frac{s}{2}}), \ \ \ s\ge 0.$$
This is a Hilbert space with inner-product $\langle A^{\frac{s}{2}}u,
A^{\frac{s}{2}} v \rangle$ for $u,v \in D(A^{\frac{s}{2}})$.  The
norm $|A^{\frac{s}{2}}u|$ is equivalent to the $H^s$ norm.

We first prove the local well-posedness of classical solutions.
\begin{thm} \label{thm1}
For $s>\frac{n}{2} +1$, and $u_0 \in W^s$, $d_0 \in H^{s+1}_h(TM)$,
there exists $T>0$ depending only on
the data and $M$, such that
$$ u \in C^0\left([0,T], W^s\right), \ \ \ d\in C^0\left([0,T], H^{s+1}_h(TM)
\right)$$
are solutions to the system of equations (\ref{lc}).
\end{thm}
\begin{proof}
It will be convenient to recast the equations (\ref{lc.c}) and (\ref{lc.e})
so that the solution has zero trace on ${\partial M}$.  For any boundary
data $h\in H^{s+{\frac{1}{2}}}(T{\partial M})$, we may choose $\psi \in 
H^{s+1}(TM)$ such that $\operatorname{trace}(\psi)=h$.  Let
$$\tilde{d} = d - \psi, \text{ so that } \tilde{d}|_{{\partial M}}=0.$$
We rewrite the system (\ref{lc}) as an evolution equation
in
$$ X^s \equiv W^s \oplus H_0^{s+1}(TM):$$
\begin{subequations}
  \label{lc2}
\begin{gather}
u_t + \nu Au + P\nabla_uu = 
-P\lambda\operatorname{Div}(\nabla [\tilde{d}+\psi]^T\cdot 
\nabla [\tilde{d}+\psi]) \,,
         \label{lc2.a}\\
\tilde{d}_t + \nabla_u\tilde{d} = \gamma\left(\hat\Delta d -
\tilde{f}(\tilde{d}) + \hat\Delta \psi \right) \,,
         \label{lc2.b}\\
u=0 \ \rm{on} \ {\partial M}, \ \ \tilde{d} = 0 \ \rm{on} \ {\partial M} \,
\rm{\text{ or }} {\partial M} = \emptyset\,,
         \label{lc2.c}\\
u(0,x)=u_0, \ \ \tilde{d}(0,x) = \tilde{d}_0(x) \equiv d_0(x) + \psi(x)\,,
         \label{lc2.d}
\end{gather}
\end{subequations}
where
$$\tilde{f}(\tilde{d}) \equiv 
{\frac{1}{\epsilon^2}} \left( |\tilde{d}+\psi|^2 -1
\right)\left(\tilde{d}+\psi\right).$$
We define the vector 
$$x \equiv (u,\tilde d) \in X^s;$$
since $H^{s-1}$-class vector fields form a Schauder ring for 
$s>{\frac{n}{2}}+1$, we may define the maps
\begin{equation}\label{phi1}
\begin{array}{c}
\phi_1: X^s \rightarrow V^{s-1}, \\
\phi_1(x) = -P\left(\nabla_uu + \operatorname{Div}(\nabla [\tilde{d}+\psi]^T
\cdot\nabla [\tilde{d}+\psi]) \right),
\end{array}
\end{equation}
and
\begin{equation}\label{phi2}
\begin{array}{c}
\phi_2: X^s \rightarrow H^{s}(TM), \\
\phi_2(x) = -\nabla_u \tilde d - \gamma\tilde{f}(\tilde{d}) + \gamma \hat 
\Delta \psi.
\end{array}
\end{equation}
Thus, the vector
$$\Phi\equiv (\phi_1, \phi_2): X^s \rightarrow V^{s-1}\times H^s(TM).$$
We are using the fact that the projector $P$ maps $H^{s-1}$ to itself.
To see this, we write 
$$P \operatorname{Div}(\nabla d^T \cdot \nabla d) =
\operatorname{Div}(\nabla d^T \cdot \nabla d) -\operatorname{grad} q,$$
where $q$ solves the Neumann problem
$$
\begin{array}{c}
\Delta q = \operatorname{div}\operatorname{Div}(\nabla d^T \cdot \nabla d),\\
g(\operatorname{grad}q,n) = g\left(
\operatorname{Div}(\nabla d^T \cdot \nabla d),n\right)
\text{ on } {\partial M}.
\end{array}
$$
Since $\operatorname{div}\operatorname{Div}(\nabla d^T \cdot \nabla d)$ is
in $H^{s-2}(M)$ and $\operatorname{Trace}
\operatorname{Div}(\nabla d^T \cdot \nabla d)$ is in 
$H^{s-\frac{3}{2}} ({\partial M})$, by elliptic regularity $q$ is in $H^s(M)$ 
so that $\operatorname{grad}q$ is in $H^{s-1}(TM)$, as desired.  One sees that 
$P \nabla_uu $ is also in $H^{s-1}(TM)$ by a similar argument.

Next, we define the semigroup
$$S(t) =
\left[\begin{array}{cccc}
        e^{-t \nu A} & 0   \\
        0   & e^{t\gamma \hat\Delta} \\
        \end{array} \right]
$$
We can now express the system (\ref{lc2}) as the integral equation
\begin{equation}\label{integral}
x_t(t,\cdot) = S(t)x_0 - \int_0^t S(t-s) \Phi(x(s)) ds = \Psi x(t,\cdot).
\end{equation}

Since $e^{-t \nu A}: W^s \rightarrow W^s$ and $e^{t\gamma \hat\Delta}:
H^{s+1}_0(TM) \rightarrow H^{s+1}_0(TM)$ are strongly continuous semigroups,
it follows that
\begin{equation}\label{strong}
S(t):X^s \rightarrow X^s \ \text{\rm is a strongly continuous semigroup for} \
t\ge 0,
\end{equation}
and that for $t>0$, $S(t):{\mathcal V}^{s-1} \times H^s(TM) \rightarrow
X^s$;  furthermore, we have the usual estimate (see, for example, 
\cite{T})
\begin{equation}\label{S_est}
\|S(t)\|_{L(V^{s-1} \times H^s(TM),X^s)} \le C t^{-\frac{1}{2}},
\ \ \ \ t\in(0,1].
\end{equation}
Using the fact that for $s>(n/2)+1$, 
$P:H^s(TM)\cap H^{s-{\frac{1}{2}}}(T{\partial M}) 
\rightarrow V^s$ is a bounded projection,
we obtain that
\begin{equation}\label{lip}
\Phi: X^s  \rightarrow V^{s-1} \times H^s(TM) \ 
\text{\rm is a locally Lipschitz map};
\end{equation}
namely,
$$
\|\phi_1(u,\tilde{d}) - \phi(v,e)\|_{s-1}
\le C_1\left( \|u-v\|_s, \|\tilde{d} - e\|_{s+1} \right)
$$
$$\|\phi_2(u,\tilde{d}) - \phi(v,e)\|_{s}
\le C_2\left( \|u-v\|_s, \|\tilde{d} - e\|_{s+1} \right)
$$
where $C_1$ and $C_2$ depend on $\|u\|_s$, $\|v\|_s$, $\|\tilde{d}\|_{s+1}$,
$\|e\|_{s+1}$, and $\|\psi\|_s$.

Fix $\alpha>0$ and set
$$ Z = \{ x \in C([0,T],X^s) \ | \ x(0) = (u_0, \tilde{d}_0), \ 
\|x(t,\cdot) - x(0) \|_{X^s} < \alpha\}.
$$
We want to choose $T$ sufficiently small so that $\Psi:Z \rightarrow Z$
is a contraction.   By (\ref{strong}), we can choose $T_1$ so that
$$\|S(t)x_0 - x_0\|_{X^s} \le \alpha/2 \ \ \ \forall t\in [0,T_1].$$
If $x \in Z$, then by (\ref{lip}) we have a bound
$$\|\Phi(x(s))\|_{{\mathcal V}^{s-1}\times H^s(TM)} \le K_1 \ \rm{for} \
s\in [0,T_1].$$
Using (\ref{S_est}), we have that
$$
\left\| \int_0^t S(t-s) \Phi(x(s)) ds\right\|_{X^s} \le C t^{\frac{1}{2}} K_1;
$$
hence, for $t\in[0,T_2]$, and with $x=(u,\tilde{d})$,
\begin{align*}
& \|\Psi(u(t),\tilde{d}(t)) - \Psi(v(t),e(t))\|_{X^s} \\
&\qquad \qquad  =
\left\| \int_0^t S(t-s) \left[\Phi(u(s),\tilde{d}(s))
-\Phi(v(s),e(s))\right] ds\right\|_{X^s}\\
&\qquad\qquad
 \le C t^{\frac{1}{2}} K \sup \|(u(s),\tilde{d}(s)) - (v(s), e(s))\|_{X^s}.
\end{align*}
Choosing $T \le T_2$ small enough so that $CT^{\frac{1}{2}} K < 1$, we
see that by the contraction mapping theorem, $\Psi$ has a unique fixed point
in $Z$, and this proves the theorem.
\end{proof}

Using the contraction mapping theorem, we can also establish the local
well-posedness for a weaker class of solutions.
\begin{thm} \label{thm1a}
Suppose $2\le \operatorname{dim}(M) \le 5$ and set
$s_0 = {\frac{n}{4}} + {\frac{1}{2}}$.  For $s\in (s_0,2)$ and
$u_0 \in W^s$, $d_0 \in H^{s+1}_h(TM)$,
there exists $T>0$ depending only on
the data and $M$, such that
$$ u \in C^0\left([0,T], W^s\right), \ \ \ d\in C^0\left([0,T], H^{s+1}_h(TM)
\right)$$
are solutions to the system of equations (\ref{lc}).
\end{thm}
\begin{proof}
We keep the same notation as in the proof of Theorem \ref{thm1}.  For 
$s\in (s_0,2)$, we have that
$$
\|S(t)\|_{L(V^0 \times H^1(TM),X^s)} \le C t^{-\gamma},
\ \ \gamma \in (0,1), \ \ t\in(0,1].
$$
Thus, it suffices to prove that for $s\in (s_0,2)$, the map
$\Phi: X^s \rightarrow V^0 \times H^1(TM)$ is locally Lipschitz.
Using Lemma 5.3 [\cite{T}, Chapter 17], we have that for $s\in(s_0,2)$
$(f,g) \mapsto f\ g: H^s \times H^s \rightarrow H^1$.  It follows that
$u \mapsto u \otimes u: H^s \rightarrow H^1$, $d \mapsto (\nabla d^T \cdot
\nabla d): H^{s+1} \rightarrow H^1$, and $(u,d) \mapsto \nabla_ud :
H^s \times H^{s+1} \rightarrow H^1$.  The fact that $d \mapsto f(d):
H^{s+1} \rightarrow H^1$ follows because $H^{s+1}$ forms a Schauder ring.
Hence, $\Phi$ is indeed locally Lipschitz, and the remainder of the proof
is identical to the one for Theorem \ref{thm1}.
\end{proof}

\section{Basic Energy Laws on Riemannian manifolds}
\label{sec_energy}

In this section, we show that the system (\ref{lc}) admits the following
energy law:
\begin{equation}\label{energylaw}
{\frac{d}{dt}}E = -
\left( \nu |\operatorname{Def} u|^2 + \lambda \gamma|\hat\Delta d
-f(d)|^2\right) - \lambda \operatorname{Trace}
\langle R(\cdot,u)d, \nabla_{\cdot} d \rangle,
\end{equation}
where $E$ is given by (\ref{energy}), and $e_i$ denotes a local orthonormal
frame.  When $M$ has zero curvature, then $E$ is a Lyapunov function for
the system (\ref{lc}), with the property that
$$E(u(t),d(t)) \le E(u_0,d_0), \ \ \ \forall t\ge 0,$$
and if $E(u(t_1),d(t_1)) = E(u(t_2),d(t_2))$ for $t_1 < t_2$, then
$(u(t),d(t))=(u^*,d^*)$ are equilibrium solutions.  Even, when the
curvature $R \neq 0$, the energy remains uniformly bounded.

This bound, in turn, then yields an a priori uniform bound for the pair 
$(u,d)$ which 
shows that all solutions  eventually enter an absorbing ball in $W^1\times 
H^2(TM)$.

The following two lemmas are standard:
\begin{lemma}\label{lemma1}
If ${\partial M} \neq \emptyset$, then for $d_0 \in H^2_h$,
$$|d(t,\cdot)|_{L^\infty} \le 1, \ \ \forall t>0.$$
\end{lemma}
\begin{proof}
We compute the pointwise inner-product of
(\ref{lc3.b}) with $d$, and use the fact that 
$g(\hat \Delta d(x), d(x)) = (1/2) \triangle( g(d(x),d(x))) - 
g(\nabla d, \nabla d)$. Hence, $\forall x \in M$, we obtain
$${\frac{1}{2}}{\frac{d}{dt}} g(d,d) 
+ {\frac{1}{2}}g( \operatorname{grad}[ g(d,d)], u)
- {\frac{1}{2}} \triangle( g(d,d)) + 
g(\nabla d, \nabla d) = -\frac{1}{\epsilon^2} [g(d,d)^2 - g(d,d)].$$
Now suppose that $\max_{t,x} g(d(t,x),d(t,x))$ occurs at $(t_0,x_0)$, 
an element of the parabolic interior;
then 
$$\frac{d}{dt}g(d(t_0,x_0),d(t_0,x_0)) =0,
\nabla d(t_0,x_0)=0,  \operatorname{grad}[g(d(t_0,x_0),d(t_0,x_0))]=0,$$
 and 
$\operatorname{Hess} g(d(t_0,x_0),d(t_0,x_0))  <0$.  This implies that
$$ - {\frac{1}{2}} \triangle( g(d(t_0,x_0),d(t_0,x_0))) >0,$$ 
but 
$$-\frac{1}{\epsilon^2}[g(d(t_0,x_0),d(t_0,x_0))^2 - 
g(d(t_0,x_0),d(t_0,x_0))] <0,$$ 
which is a contradiction.
\end{proof}

\begin{lemma}\label{lemma2}
If ${\partial M}= \emptyset$, then for $d_0 \in H^2(M,B^\delta)$, 
$$|d(t,\cdot)|_{L^\infty} \le \delta, \ \ \forall t>0.$$
\end{lemma}
\begin{proof}
This again follows from the maximum principle above.
\end{proof}

\begin{prop}\label{prop1}
The energy law (\ref{energylaw}) holds, and there exists an absorbing set
for $(u,d) \in W^0 \times H^1_h(TM)$ if ${\partial M}\neq \emptyset$ and for 
$(u,d) \in W^0 \times H^1(M,B^\delta)$ if ${\partial M} = \emptyset$.
\end{prop}
\begin{proof}
Using the formula (\ref{R}), we rewrite (\ref{lc.a}) and (\ref{lc.c}) as
\begin{subequations} \label{lc3}
\begin{gather}
u_t + \nabla_uu = -\operatorname{grad}p + \nu \operatorname{Div}
\operatorname{Def}u - \lambda \nabla d^T \cdot \hat\Delta d
- g(R(e_i,\cdot) d, \nabla_{e_i} d) \,, \label{lc3.a}\\
d_t + \nabla_ud = \gamma\left(\hat\Delta d - 
\frac{1}{\epsilon^2}(|d|^2 -1)d\right) \,, \label{lc3.b}
\end{gather}
\end{subequations}
where $e_i$ is any local orthonormal frame.  Adding the $L^2$ inner-product
of (\ref{lc3.a}) with $u$ to the $L^2$ inner-product of (\ref{lc3.b}), we
obtain the basic energy law
\begin{align}
&{\frac{1}{2}}{\frac{d}{dt}}\left( |u|^2 + \lambda |\nabla d|^2
+ 2 \lambda \int_M F(d(x)) \mu \right) \nonumber\\
&\qquad \qquad
= -\left( \nu |\operatorname{Def} u|^2 + \lambda \gamma|\hat\Delta d
-f(d)|^2\right)
- \lambda \operatorname{Trace}
\langle R(\cdot,u)d, \nabla_{\cdot} d \rangle.
\label{energy1}
\end{align}
In the case of a flat manifold, such as a bounded domain in ${\mathbb R}^n$,
$R=0$, and (\ref{energy1}) reduces to the basic energy law (1.8) in
\cite{LL1}.  

From Lemmas \ref{lemma1} and \ref{lemma2}, we have that 
\begin{equation}\label{max}
|d(t,\cdot)|_{L^\infty} \le  C, \ \ \ \ t > 0.
\end{equation}
It follows that
\begin{align*}
\operatorname{Trace} \langle R(\cdot,u)d, \nabla_{\cdot} d \rangle
&\le C |R|_{L^\infty} |u| |d|^{\frac{1}{2}} |\hat\Delta d|^{\frac{1}{2}} \\
&\le C \varepsilon |\hat\Delta d|^2 + {\frac{C}{\varepsilon}}
\left( |M| |R|_{L^\infty} |u| \right)^{\frac{4}{3}} \\
&\le C \varepsilon |\hat\Delta d|^2 +  \varepsilon |u|^2 +  
{\frac{C}{\varepsilon^4}} \left( |M| |R|_{L^\infty} \right)^4 \\
&\le C \varepsilon |\hat\Delta d|^2 +  c_0^{-1}(M) \varepsilon 
|\operatorname{Def}u|^2 +  
{\frac{C}{\varepsilon^4}} \left( |M| |R|_{L^\infty} \right)^4\, ,
\end{align*}
where the second and third inequalities follow from Young's inequality
(\ref{young}), and the last inequality follows from the Poincar\'{e}
inequality for $c_0(M)>0$, a positive constant depending on $M$.
Taking $\varepsilon>0$ sufficiently small so that
$$K= \min (c_0 - \varepsilon, 1- 2 \varepsilon) >0,$$
the basic energy law (\ref{energy1}) on a Riemannian manifold yields the
following differential inequalities:

\begin{subequations} \label{I}
\begin{align}
& \frac{1}{2}\frac{d}{dt} \left[ |u|^2 + |\nabla d|^2 + 2\int_M F(d)\mu\right]
\nonumber \\
&\qquad \qquad \le -K \ C \left[
 |\operatorname{Def}u|^2 + |\hat\Delta d|^2 + 2\int_M F(d)\mu\right]+\rho_0 \,,
\label{I.a}\\
& \frac{1}{2}\frac{d}{dt} \left[ |u|^2 + |\nabla d|^2 + 2\int_M F(d)\mu\right]
\nonumber \\
&\qquad \qquad \le -K \ C \left[
 |u|^2 + |\nabla d|^2 + 2\int_M F(d)\mu\right] +\rho_0 \,,
\label{I.b}
\end{align}
\end{subequations}
where
$$\rho_0 = C\left[( K + C/\varepsilon -1)|M| + \frac{1}{\varepsilon^4}|M|^4
|R|_{L^\infty}\right].$$
Using the classical Gronwall lemma, we obtain 
\begin{align*}
&\left[ |u|^2 + |\nabla d|^2 + 2\int_M F(d)\mu\right] \\
&\qquad \qquad \le
\left[|u_0|^2 + |\nabla d_0|^2 + 2\int_M F(d_0)\mu\right] e^{-K\, Ct}
+\rho_0(1-e^{-K\, Ct}).
\end{align*}
Thus,
\begin{subequations} \label{star}
\begin{gather}
\limsup_{t \rightarrow \infty}
\left[ |u(t)|^2 + |\nabla d(t)|^2 + 2\int_M F(d(t))\mu\right]  \le \rho_0
\label{star.a}\\
\limsup_{t \rightarrow \infty}
\left[ |u(t)|^2 + |\nabla d(t)|^2 \right]  \le \rho_0+ 2|M|.
\label{star.b}
\end{gather}
\end{subequations} 
\end{proof}

When $R=0$, we do not need to rely on the maximum principle to establish
Proposition \ref{prop1} or to establish the existence of an $L^\infty$ absorbing
set for (\ref{lc3.b}).  

\begin{lemma}\label{lemma3}
If $M=\emptyset$ and $R=0$, then for $d_0 \in H^2(TM)$, there exists
$\rho_\infty >0$ and some $t^* >0$ independent of $d_0$ such that
$$ |d(t,\cdot)|_{L^\infty} \le \rho_\infty \ \ \forall t > t^*.$$
\end{lemma}
\begin{proof}
When $R=0$, from the energy law (\ref{energylaw}), we see that
there exists and $L^2$ absorbing set, so that for some $t > t^*$  all bounded
subsets of $L^2(TM)$ will enter the $L^2$ ball of radius $\rho_0$.

For $p>2$,
we take the pointwise inner-product of
(\ref{lc3.b}) with
$p|d|^{p-2}d$ and integrate over $M$ to obtain the
differential inequality
\begin{align}
\frac{d}{dt}|d|_{L^p}^p &= -p\int_M |\nabla d|^2 |d|^{p-2} \mu
- p(p-2)\int_M |d|^{p-2} | \nabla|d| |^2 \mu \nonumber\\
&\qquad \qquad+{\frac{1}{\epsilon^2}} \left(|d|_{L^p}^p - |d|_{L^p}^{p-2}\right)
\nonumber\\
& \le
- p(p-2)\int_M |d|^{p-2} | \nabla|d| |^2 \mu
+\frac{1}{\epsilon^2} |d|_{L^p}^p. \label{diffeq}
\end{align}
Using the interpolation inequality (see \cite{LS} for details and further
applications)
$$ |d|_{L^p}^p \le C_p |d|^2 \left( \int_M |\nabla |d|^{\frac{p}{2}}|^2 \mu
\right)^{\frac{p-2}{p}} =
C_p|d|^2 \left( {\frac{p^2}{4}} \int_M |d|^{p-2} |\nabla |d||^2 \mu
\right)^{\frac{p-2}{p}},
$$
we see that
$$-p(p-2) \int_M |d|^{p-2} | \nabla|d| |^2 \mu
\le \left[ \left( \frac{1}{C_p}\right)^{\frac{p}{2(p-2)}}
\frac{4p(p-2)}{p^2} \rho_0^{\frac{p}{p-2}}\right] \left( |d|_{L^p}^p \right)
^{\frac{p}{p-2}}.$$
Using Bernoulli's trick in the differential inequality (\ref{diffeq}), we
get a uniform bound for $|d(t,\cdot)|_{L^p}$ which is independent of $p$ (even
if the constant $C_p$ tends to infinity), and
thus we may pass to the limit as $p \rightarrow \infty$.
\end{proof}
Using Lemma \ref{lemma3} we immediately have

\begin{prop}\label{prop1a}
If $R=0$ and ${\partial M} = \emptyset$, we have the energy law
$$
{\frac{d}{dt}}E = -
\left( \nu |\operatorname{Def} u|^2 + \lambda \gamma|\hat\Delta d
-f(d)|^2\right),
$$
and there exists an absorbing set
for $(u,d) \in W^0 \times H^1(TM)$.
\end{prop}

\begin{prop} \label{H1H2}
For $\operatorname{dim}(M)=2$, there exists an absorbing set for
$(u,d) \in W^1 \times H^2_h(TM)$ if ${\partial M} \neq \emptyset$,
for $(u,d) \in W^1 \times H^2(M,B^\delta)$ if ${\partial M} =\emptyset$ and
$R\neq 0$, and for $(u,d) \in W^1 \times H^2(TM)$ if ${\partial M}=\emptyset$
and $R=0$.
\end{prop}
\begin{proof}
It follows from (\ref{I.a}) that
\begin{align*}
& K\, C \int_t^{t+r}
\left[ |\operatorname{Def} u(s)|^2 + |\hat \Delta d(s)|^2 + 
2\int_M F(d(s))\mu\right] ds \\
& \qquad \qquad \qquad \qquad \le
r \rho_0 +
\left[ |u|^2 + |\nabla d|^2 + 2\int_M F(d))\mu\right], \ \ \ \forall r>0,
\end{align*}
so 
$$
\limsup_{t \rightarrow \infty}\int_t^{t+r}
\left[ |\operatorname{Def} u(s)|^2 + |\hat \Delta d(s)|^2 + 
2\int_M F(d(s))\mu\right] ds  \le (r+1) \rho_0.
$$
Therefore,
\begin{equation}\label{9}
\int_t^{t+r}
\left[ \left|\operatorname{Def} u(s)\right|^2 + |\hat \Delta d(s)|^2 
+ 
2\int_M F(d(s))\mu\right] ds \text{ is uniformly bounded}.
\end{equation}

Now let
$$A^2 = |\operatorname{Def} u|^2 + |\hat \Delta d -f(d)|^2, \ \
B^2 = |\nabla \operatorname{Def} u|^2 + | \nabla 
(\hat \Delta d -f(d))|^2.
$$
Using (\ref{9}) and (\ref{max}), we have that
\begin{equation}\label{A}
\int_t^{t+r} A^2(s) ds \ \ \text{ is uniformly bounded.}
\end{equation}
In the case that $R=0$, it follows from a similar argument as in 
(4.4)-(4.8) of \cite{LL1} that
for some constants $c_1, c_2, c_3 >0$,
\begin{equation}\label{LL}
\frac{d}{dt} A^2(t) + c_1 B^2(t) \le  c_2 A^4(t) +c_3.
\end{equation}
When $R \ne 0$, we find that  for $c_4>0$,
$$
\frac{d}{dt} A^2(t) + c_1 B^2(t) \le  c_2 A^4(t) +c_3 + c_4
\operatorname{Trace}\langle R( \cdot, \triangle u)
d, \nabla_\cdot d \rangle.$$
The last term is bounded by
$\varepsilon |R|_{L^\infty}^2 |\operatorname{Div}\operatorname{Def} u|^2
+(C/\varepsilon)|\nabla d|^2$, so
by taking $\varepsilon >0$ sufficiently small and adjusting the
constants as necessary, we see that (\ref{LL}) still holds.

Thus, using (\ref{star}) and appealing to the uniform Gronwall lemma
(see, for example, \cite{Tem}), we see that
$$
A(t) \ \ \ \text{is uniformly bounded in time.}
$$
Because of (\ref{max}), we may extract a uniform bound for
$|\operatorname{Def} u|^2 + |\hat\Delta d|^2$.  Hence, we have an a priori
uniform bound for $u$ in the $H^1$ topology and for $d$ in the $H^2$
topology.
\end{proof}

\section{Global Well-posedness and Global Attractors}
\label{sec_global}

We shall first consider a closed Riemannian manifold such as, for example, the
two-sphere ${\mathbb S}^2$; for such manifolds, simple brute-force energy 
estimates work.
\begin{thm}\label{thm2}
For $n=2$, $s>1$, ${\partial M}= \emptyset$, and $u_0 \in W^s$, $d_0 \in 
H^{s+1}(M,B^\delta)$, 
$$ u \in C^0([0,\infty], W^s), \ \ \ d\in C^0([0,\infty], H^{s+1}(M,B^\delta))$$
are solutions to the system of equations (\ref{lc}).  Moreover, there
exists a compact global attractor for the system (\ref{lc}) 
in $W^{s-1} \times H^{s}(M,B^\delta)$.
In the case that $R=0$, we can replace $H^{s+1}(M,B^\delta)$ with 
$H^{s+1}(TM)$.
\end{thm}
\begin{proof}
Taking the  $H^s$ inner-product of (\ref{lc.a}) with $u$ and  adding
the $H^{s+1}$ inner-product of (\ref{lc.c}) with $d$, we find that
\begin{align}
{\frac{1}{2}}\frac{d}{dt}\left(|u|^2_s + |d|^2_{s+1}\right)
& \le 
-\nu |u|^2_{s+1} - \gamma |d|^2_{s+2} 
+ \lambda \left|\left\langle P \operatorname{Div}(\nabla d^T \cdot \nabla d), 
u\right\rangle_s \right|
\nonumber\\
& \qquad
+ \left|\left\langle P \nabla_uu, u \right\rangle_s\right| 
+ \left|\left\langle \nabla_ud, d \right\rangle_{s+1}\right|
+ \gamma\left|\left\langle f(d), d \right\rangle_{s+1} \right|.
\label{basic}
\end{align}
We shall estimate each of the nonlinear terms on the right-hand-side of 
(\ref{basic});  as we showed in the proof of Theorem \ref{thm1}, the projection
$P$ acting on the nonlinear terms, maps $H^{s-1}$ into itself continuously, so
it suffices to estimate $\langle \nabla_uu, u \rangle_s$ and
$\langle\operatorname{Div}(\nabla d^T \cdot \nabla d), u \rangle_s$ in the
third and fourth terms.  Using Proposition \ref{H1H2}, we may interpolate the
nonlinear terms between $|u|_1$ and $|u|_{s+1}$ and $|d|_2$ and $|d|_{s+2}$,
respectively.

We have that
\begin{align*}
\langle \nabla_uu, u \rangle_s & = \sum_{\alpha=s} \langle D^\alpha(\nabla_uu),
D^\alpha u \rangle = 
\sum_{\alpha=s} \sum_{\stackrel{|\beta| \le s }{\alpha -\beta\ge 0 }}
\langle c_{\alpha,\beta} (D^\beta \nabla u) (D^{\alpha-\beta}u), D^\alpha u 
\rangle \\
& =
\sum_{\alpha=s} \sum_{\stackrel{|\beta| \le {s-1} }{\alpha -\beta\ge 0 }}
\langle c_{\alpha,\beta} (D^\beta \nabla u) (D^{\alpha-\beta}u), D^\alpha u 
\rangle \\
& \le  C
\sum_{m=0}^{s-1} |D^{m+1}u|_{L^4} |D^{s-m}|_{L^4}|D^su|.
\end{align*}
where we set $|\beta|=m$ so that $|\alpha - \beta|=s-m$,  and the 
last equality follows from the fact that $\langle \nabla_u(D^su),
D^su\rangle =0$, since $\operatorname{div}u=0$.  For
$m=0,...,s-1$, we use (\ref{interp2}) and (\ref{interp3}) to estimate

\begin{align*}
|D^{m+1}u|_{L^4} |D^{s-m}|_{L^4}|D^su| & \le
C |D^{m+1}u|^{\frac{1}{2}}
|D^{m+2}u|^{\frac{1}{2}}
|D^{s-m}u|^{\frac{1}{2}}
|D^{s-m+1}u|^{\frac{1}{2}} |D^s u| \\
\le C |u|_1^{\frac{s+1}{s}}
|u|_{s+1}^{\frac{2(s-1)}{2s}}.
\end{align*}
Using Young's inequality, 
\begin{equation}\label{young}
a^\lambda b \le \varepsilon a + {\frac{C}{\varepsilon}} b^{\frac{1}{1-\lambda}},
\ \ \ \ a,b>0,\ \ \ \ 0 < \lambda < 1,
\end{equation}
it follows that
$$ \langle P \nabla_uu, u \rangle_s \le  \varepsilon |u|^2_{s+1} +
{\frac{C}{\varepsilon}} |u|_1^{s+1}.$$

For the next term, we have that
\begin{align}
\langle \operatorname{Div}(\nabla d^T \cdot \nabla d), u \rangle_s & =
\sum_{|\alpha|=s} \langle D^\alpha \operatorname{Div}(\nabla d^T \cdot
\nabla d), D^\alpha u \rangle \nonumber \\
&\le 
C \sum_{\stackrel{|\alpha|=s}{|\gamma|=s+1}} 
\sum_{\stackrel{|\beta|\le s}{\gamma-\beta\ge 0}}
\int_M (D^\beta \nabla d) (D^{\gamma-\beta}\nabla d) (D^\alpha u) \mu 
\nonumber \\
&\le C \sum_{m=0}^{s+1} \langle(D^{m+1}d) (D^{s-m+2}d), (D^s u)\rangle, 
\label{2}
\end{align}
where $m= |\beta|$.  In the case that $m=1,...,s$,  (\ref{2}) is bounded by
\begin{align}
C\sum_{m=1}^s |D^{m+1}d|^2_{L^4} |D^{s-m+1}d|^2_{L^4} 
&\le C |d|_2 |d_{s+2} |u|_1^{\frac{1}{s}} |u|_{s+1}^{\frac{s-1}{s}}\nonumber\\
&\le \varepsilon |d|_{s+2}^2 + {\frac{C}{\varepsilon}} |d|_2^2 |u|_s
^{\frac{2}{s}}|u|_{s+1}^{\frac{2(s-1)}{s}}, \label{3}
\end{align}
where the first inequality follows from repeated use of (\ref{interp3}),
and the last inequality follows from $ab \le \varepsilon a^2 + 
(C/\varepsilon)b^2$,
where $a,b>0$.  One more application of (\ref{young}) shows that
(\ref{3}) is bounded by
$$ \varepsilon|d|_{s+2}^2 + \varepsilon |u|_{s+1}^2 + 
{\frac{C}{\varepsilon^{1+s}}} |d|_2^{2s} |u|_1^2.$$
In the case that $m=0,s+1$, (\ref{2}) is bounded by
\begin{align*}
C|d|_{s+2} |d|_{L^4} |D^s u|_{L^4} & \le
C |d|_{s+2} |d|_2 |u|_1^{\frac{1}{2s}} |u|_{s+1}^{\frac{2s-1}{2s}} \\
&\le \varepsilon |d|_{s+2}^2 + {\frac{C}{\varepsilon}} |d|_2^2|u|_1^
{\frac{1}{s}} |u|_{s+1}^{\frac{2(s-1)}{s}} \\
&\le \varepsilon |d|_{s+2}^2 + \varepsilon |u|_{s+1}^2 +
{\frac{C}{\varepsilon^{1+s}}} |d|_2^{2s} |u|_1,
\end{align*}
where the first inequality follows from (\ref{interp3}), and the
last two inequalities follow from (\ref{young}).  It follows that
$$ \langle P \operatorname{Div}(\nabla d^T \cdot \nabla d), u \rangle_s 
\le \varepsilon |u|^2_{s+1} + \varepsilon|d|^2_{s+2} 
+{\frac{C}{\varepsilon^{1+s}}} |d|_2^{2s}\left( |u|_1 + |u|_1^2\right).$$
 
We next compute that
\begin{align}
\langle \nabla_u d, d \rangle_{s+1} & =
\sum_{|\alpha|={s+1}} \langle D^\alpha (\nabla_u d), D^\alpha u \rangle 
\nonumber \\
&\le 
C \sum_{|\alpha|={s+1}} \sum_{\stackrel{|\beta|\le {s+1}}{\alpha-\beta\ge 0}}
\langle (D^\beta \nabla d) (D^{\alpha-\beta}u), D^\alpha \nabla d
\rangle \nonumber\\
&\le C\sum_{m=0}^s \langle (D^{m+1}d)(D^{s-m+1}u), D^{s+1}d \rangle,
\label{1}
\end{align}
since for $m=s+1$, we have that $\langle \nabla_u (D^{s+1}d), D^{s+1}d\rangle
=0$.  We estimate the case $m=0$ first in (\ref{1}):
\begin{align*}
|Dd|_{L^4} |D^{s+1}d|_{L^4} |D^{s+1}u| & \le
C |d|_2 |u|_{s+1} |d|_{s+1}^{\frac{1}{2}}
|d|_{s+2}^{\frac{1}{2}} \\
&\le C |d|_2^{\frac{2s+1}{2s}} |u|_{s+1} |d|_{s+2}^{\frac{2s-1}{2s}} \\
&\le \varepsilon |d|^2_{s+2} + {\frac{C}{\varepsilon}}
\left(|d|_2^{\frac{2s+1}{2}} |u|_{s+1}\right)^{\frac{4s}{2s+1}} \\
&\le \varepsilon |d|^2_{s+2} + \varepsilon |u|_{s+1}^2 +
{\frac{C}{\varepsilon^{2s+2}}}|d|_2^{4s+2},
\end{align*}
where the last two inequalities follow from two applications of
the Young's inequality.

For the cases $1\le m \le s$, 
(\ref{1}) is bounded by $C |D^{m+1}d|_{L^4}
|D^{s-m+1}u|_{L^4} |d|_{s+1}$, so by (\ref{interp2})  and (\ref{interp3}), we 
find that for $m=1,...,s$,
\begin{align*}
&|D^{m+1}d|^{\frac{1}{2}} |D^{m+2}d|^{\frac{1}{2}} |D^{s-m+1}u|^{\frac{1}{2}}
|D^{s-m+2}u|^{\frac{1}{2}} |D^{s+1}d| \\
& \qquad \qquad \qquad \qquad 
\le C |d|_2^{\frac{2s-2m+3}{2s}} |d|_{s+2}^{\frac{2m+2s-3}{2s}}
|u|_1^{\frac{2m-1}{2s}} |u|_{s+1}^{\frac{2s-2m+1}{2s}} \\
& \qquad \qquad \qquad \qquad 
\le \varepsilon |d|_{s+2}^2 + \left(\frac{C}{\varepsilon}
|d|^2_2 |u|_1^{\frac{4m-2}{2s-2m+3}}\right)
\left( |u|_{s+1}^{\frac{4s-4m+2}{2s-2m+3}}\right),
\end{align*}
where we have used Young's inequality for the last step.  Another application
of Young's inequality yields the estimate
$$ \langle \nabla_ud, d \rangle_{s+1} \le 
\varepsilon |d|_{s+2}^2 + \varepsilon |u|_{s+1}^2 +
{\frac{C}{\varepsilon^{{\frac{2s-2m+5}{2s-2m+3}}}}} |d|_2^{\frac{4}{2s-2m+3}}
|u|_1^{2m-1}, \ \ \ m=1,...,s.$$

For the final nonlinear term, we have that
\begin{align*}
\langle f(d), d\rangle_{s+1} & \le  C \sum_{m=0}^{s+1}
\sum_{n=0}^{m} \langle (D^n d)(D^{m-n}d) (D^{s+1-m}d), D^{s+1}d \rangle \\
&\le |D^n d|_{L^8} |D^{m-n}d|_{L^8} |D^{s+1-m}d|_{L^4} |D^{s+1}d|.
\end{align*}
Using the estimate
\begin{equation}\label{interp4}
|v|_{L^8} \le |v|_{L^2}^{\frac{5}{8}} |v|_2^{\frac{3}{8}},
\end{equation}
together with (\ref{interp3}) and Young's inequality, we have that
$$
\langle f(d), d\rangle_{s+1}  \le  C |d|_2^{\frac{2(s+2)}{s}}
|d|_{s+2}^{\frac{2(s-2)}{s}}
\le \varepsilon |d|^2_{s+2} + {\frac{C}{\varepsilon}}|d|_2^{s+2} .$$

Letting
\begin{align*}
&\rho=C \left[\frac{1}{\varepsilon} \left(|d|_2^{s+2}+ |u|_1^{s+1} \right) 
+
\frac{1}{\varepsilon^{1+s}} |d|_2^{2s}\left(|u|_1+|u|_1^2\right)
+
\frac{1}{\varepsilon^{2+2s}} (|d|_2^{4s+2}) \right.\\
& \qquad\qquad\qquad\qquad\qquad\qquad\qquad \left.
+
\sum_{m=1}^s
\frac{1}{\varepsilon^{\frac{2s-2m+5}{2s-2m+3}}} |d|_2^{\frac{4}{2s-2m+3}}
|u|_1^{2m-1}\right],
\end{align*}
and taking $\varepsilon >0$ sufficiently small so that
$$ K = \min\left(\nu - 4 \varepsilon, \gamma - 4 \varepsilon\right) >0,$$
the basic inequality (\ref{basic}) takes the form
\begin{align*}
\frac{d}{dt}\left(|u|^2_s + |d|^2_{s+1} \right) 
& \le -K\left(|u|^2_{s+1} + |d|^2_{s+2} \right) + C\,\rho \\
& \le -K\left(|u|^2_{s} + |d|^2_{s+1} \right) + C\,\rho.
\end{align*}
Letting $c_1=CK>0$, the classical Gronwall lemma gives
$$
\left(|u|^2_s + |d|^2_{s+1} \right)  \le
\left(|u_0|^2_s + |d_0|^2_{s+1} \right) e^{-c_1t} +C\, \rho(1-e^{-c_1t}),
$$
so that
\begin{equation}\label{bound}
\limsup_{t \rightarrow \infty}
\left(|u|^2_s + |d|^2_{s+1} \right) \le C\, \rho.
\end{equation}

Thus, since the time interval of existence from Theorem \ref{thm1} only depends
on the initial data, the a priori bound (\ref{bound}) together with the
continuation property gives the global well-posedness result.  
Moreover, because of the absorbing sets that exist in $W^s \times
H^{s+1}(TM)$ by virtue of
(\ref{bound}), we obtain using
Theorem I.1.1 of \cite{Tem}, the global attractor that we asserted.
\end{proof}

For a Riemannian manifold with boundary, the above (brute force) 
$H^s$ energy
estimate does not work, because boundary terms 
arising from integration by parts 
on the diffusion term $\nu\operatorname{Div} \operatorname{Def} u$ do not
vanish.  It is possible, however, to obtain estimates on
$u_t$ and $d_t$ which provide the global well-posedness result.

We have that
\begin{align*}
u_{tt} &= 
\nu\operatorname{Div}\operatorname{Def} u_t - \nabla_{u_t}u - \nabla_uu_t 
- \operatorname{grad}p_t -\lambda \nabla d_t^T\cdot \hat\Delta d
- \lambda \nabla d^T \cdot\hat \Delta d_t \\
& \qquad - \lambda g(R(e_i,\cdot)d_t, \nabla _{e_i} d)
- \lambda g(R(e_i,\cdot)d, \nabla _{e_i} d_t)
\end{align*}
and
$$
d_{tt} = -\nabla_{u_t}d - \nabla_ud_t + \gamma \hat\Delta d_t - \gamma
\operatorname{grad}f(d) \cdot d_t.
$$
Since $d_t=0$ on ${\partial M}$ we see that 
${\frac{1}{2}}\frac{d}{dt}\left(|u_t|^2 + |\nabla d_t|^2\right) =
\langle u_t, u_{tt}\rangle - \langle \hat \Delta d_t, d_{tt}\rangle$.
Standard interpolation combined with Young's inequality yields, for constants
$c_1,c_2>0$,
$$
{\frac{1}{2}}\frac{d}{dt}\left(|u_t|^2 + |\nabla d_t|^2\right) =
c_1\left( |\operatorname{Def} u|^2 + |\hat\Delta d|^2 + c_2\right)
\left( |u_t|^2 + |\nabla d_t|^2 \right).
$$
By Proposition \ref{H1H2}, for each $t$, 
\begin{equation}\label{s1}
u(t,\cdot) \in W^1 \text{ and } d(t,\cdot) \in H^2(TM).
\end{equation}
It follows that  if $u_0 \in W^2$,  $d_0 \in H^3(TM)$, and $h\in 
H^{\frac{5}{2}}(T{\partial M})$, then $u_t(0) \in L^2(TM)$ and 
$d_t(0) \in H^1_0(TM)$ so that 
\begin{equation}\label{s2}
u_t\in L^\infty((0,\infty), W^0) \text{ and } d_t\in L^\infty(
(0,\infty),H^1_0(TM)).
\end{equation}
From (\ref{s1}), we claim that $\nabla_ud$ is in $H^\delta(TM)$ for 
$\delta \in (0,{\frac{5}{16}})$.  To see this, note that for $\varepsilon >0$
$$ 
w \mapsto w \ w : H^p \rightarrow H^{\theta(1+\varepsilon)}, \text{ where }
p = {\frac{1}{2}} + {\frac{1}{2}}\theta(1+\varepsilon) + \varepsilon \theta.
$$
We set $\delta = \theta(1+\varepsilon)$, and, for example, set
$\varepsilon ={\frac{1}{4}}$ and $\theta\le{\frac{1}{4}}$;  then $\delta\in (0,
{\frac{5}{16}})$ and $p\ge {\frac{23}{32}}$,  
so the claim is established.  Using standard elliptic regularity on
equation (\ref{lc3.b}), we see that $d\in H^{2+\delta}(TM)$, and
the $H^{2+\delta}$-norm of $d$ only depends on the initial data and $M$.
This shows that $\operatorname{Div}(\nabla d^T \cdot \nabla d)$ is in
$L^2$ so that with (\ref{s2}), we see that $u$ is in $W^2$.  By
bootstrapping, we find that $d$ is in $H^3$,  and the continuation argument
shows that the unique solution may be continued for all time.
If $h \in C^\infty (T{\partial M})$, then both $u$ and $d$ are in
$C^\infty ( (0,\infty) \times M)$.  Proposition \ref{H1H2} together with
Theorem I.1.1 of \cite{Tem} proves the existence of the global attractor
in $W^0 \times H^1_h(TM)$. Thus, we have the following
\begin{thm}\label{thm2a}
Suppose that $u_0 \in W^2$ and $d_0 \in H^3_h(TM)$.
Then there exists a unique solution
$$ u \in L^\infty( (0,\infty), W^2) \text{ and }
 d \in L^\infty( (0,\infty), H^3_h(TM)).
$$
If $h \in C^\infty(T{\partial M})$, then both $u$ and $d$ are in
$C^\infty( (0,\infty) \times M)$.  Furthermore, there exists a compact global
attractor in $W^0 \times H^1_h(TM)$.
\end{thm}

\section{Lagrangian averaged liquid crystals}
\label{sec_averaged}

As we described in the introduction, the director field $d$ describes locally 
the 
averaged direction of the constituent molecules; it is thus reasonable,
and of practical and computational importance, to locally 
average the Navier-Stokes fluid motion as well.  Recently, the Lagrangian
averaged Navier-Stokes (LANS) equations were introduced as a model
for the large scale Navier-Stokes  fluid motion which averages or filters
over the small, computationally unresolvable spatial scales (see \cite{MS1} and
the references therein).  The LANS equations are parameterized by a small
spatial scale $\alpha>0$ -- fluid motion at spatial scales smaller than 
$\alpha$ is averaged or filtered-out.    There are two types of Lagrangian
averaged Navier-Stokes equations: the {\it isotropic} and the {\it anisotropic}
versions.  We shall begin with the isotropic theory, and for simplicity of
presentation, we shall assume that $M$ is flat.

The isotropic LANS equations for the {\it mean velocity} $u(t,x)$ are given
by
\begin{subequations}
  \label{LANS1}
\begin{gather}
\partial_t(1-\alpha^2\Delta)u + \nabla_u(1-\alpha^2\Delta)u - \alpha^2
\nabla u^T \cdot \triangle u = -\operatorname{grad} p + \nu(1-\alpha2\Delta)
\Delta u
         \label{LANS1.a}\\
 \operatorname{div}u(t,x)=0 \,,
         \label{LANS1.b}\\
u=0 \ \rm{on} \ {\partial M}\,,
         \label{LANS1.c}\\
u(0,x)=u_0.
         \label{LANS1.d}
\end{gather}
\end{subequations}
Equation (\ref{LANS1.a}) has an equivalent representation as
\begin{subequations}
  \label{LANS2}
\begin{gather}
\partial_tu + \nabla_u u + {\mathcal U}^\alpha(u) = -(1-\alpha^2\Delta)^{-1}
\operatorname{grad}p + \nu \Delta u
         \label{LANS2.a}\\
{\mathcal U}^\alpha(u) = \alpha^2(1-\alpha^2 \Delta)^{-1} \operatorname{Div}
\left[ \nabla u \cdot \nabla u^T + \nabla u \cdot \nabla u - \nabla u^T \cdot
\nabla u\right].
         \label{LANS2.b}
\end{gather}
\end{subequations}
When ${\partial M}=\emptyset$ the LANS equations take on a particularly 
familiar ``sub-grid-stress'' form with (\ref{LANS2}) becoming
\begin{subequations}
  \label{LANS3}
\begin{gather}
\partial_tu + \nabla_u u + \operatorname{Div}\tau^\alpha(u) = 
- \operatorname{grad}p + \nu \Delta u
         \label{LANS3.a}\\
\tau^\alpha(u) = \alpha^2(1-\alpha^2 \Delta)^{-1} 
\left[ \nabla u \cdot \nabla u^T + \nabla u \cdot \nabla u - \nabla u^T \cdot
\nabla u\right],
         \label{LANS3.b}
\end{gather}
\end{subequations}
where $\tau^\alpha$ representing the sub-grid or ``Reynolds stress.''

The remarkable feature of the LANS equations is that, unlike the Reynolds
averaged Navier-Stokes (RANS) equations or Large Eddy Simulation (LES) models
of turbulence, no additional dissipation is put into the system.  In
fact, when $\nu =0$, the LANS equations conserve the Hamiltonian structure of
the Euler equations with  both  a modified {\it kinetic energy}
\begin{equation}\label{LANS-energy}
E^\alpha = 
{\frac{1}{2}} \int_M \left( |u|^2 + 2 \alpha^2|\operatorname{Def} u|^2
\right) \mu
\end{equation}
and {\it helicity}
\begin{equation}\label{helicity}
H^\alpha=
\int_M w \wedge dw, \ \ \ \ \ w = (1-\alpha^2\Delta) u^\flat, \ \
u^\flat = g(u, \cdot),
\end{equation}
being conserved.

This is easiest to see from equation (\ref{LANS3}), where the only term that
is added (to the NS equations) is $\operatorname{Div}\tau^\alpha(u)$;
it is precisely this term which averages the small scales, and this is
accomplished by the use of {\it nonlinear  dispersion} as opposed to 
dissipation.  A simple computation, which requires taking the $L^2$ 
inner-product of the LANS equations with $u$ when $\nu$ is set to zero, shows
that (\ref{LANS-energy}) is conserved.  Why is it so important not to
over-dissipate the NS equations?  The answer is twofold:  first, the
addition of artificial dissipation obviously and spuriously removes crucial
small-scale features, and second, artificial viscosity, which is present
in RANS or LES models, suppresses intermittency, a fundamental feature of fluid
turbulence.

Mathematically, for all $\alpha >0$, the three-dimensional LANS equations are 
globally well-posed (see \cite{MS2}), yet when the averaging parameter
$\alpha$ is taken sufficiently small, computational simulations of 
LANS are statistically
indistinguishable from the simulations of the NS equations.  Furthermore,
the LANS equations provide a tremendous computational savings as shown in
simulations of both forced and decaying turbulence (\cite{Chen}, \cite{CTR}).
Finally, the LANS equations arise from a variational principle in the
same fashion as the NS equations.  We shall therefore base our development
of the averaged liquid crystal equations on the LANS model, and introduce
the following system of equations:
\begin{subequations}
  \label{avg_lc}
\begin{gather}
u_t  -\nu\Delta u +\nabla_uu -{\mathcal U}^\alpha(u)= -(1-\alpha^2\Delta)^{-1}
\left[ \operatorname{grad}p + \operatorname{Div}(\nabla d^T \cdot \nabla d)
\right] \,,
         \label{avg_lc.a}\\
 \operatorname{div}u(t,x)=0 \,,
         \label{avg_lc.b}\\
d_t + \nabla_ud = \gamma\left(\Delta d -
\frac{1}{\epsilon^2}(|d|^2 -1)d\right) \,,
         \label{avg_lc.c}\\
u=0 \ \rm{on} \ {\partial M}, \ \ d = h \ \rm{on} \ {\partial M} \ \ \
g(h,h)=1 \, \rm{\text{ or }} {\partial M} = \emptyset\,,
         \label{avg_lc.d}\\
u(0,x)=u_0, \ \ d(0,x) = d_0 \text{ and } d_0|_{\partial M}=h \text{ if }
   {\partial M} \neq \emptyset \,,
         \label{avg_lc.e}
\end{gather}
\end{subequations}
where $\Delta$ denotes the componentwise Laplacian.

\medskip

\noindent
{\bf Averaged Variational Principle.}  Following the notation of Section
\ref{sec_variational}, we define the averaged action function
$S^\alpha: {\mathcal D}_{\mu,D}^s \times H^{s+1}(TM) \cap H^1_0(TM) 
\rightarrow {\mathbb R}$ by
\begin{align}
S^\alpha(\eta, d) = & {\frac{1}{2}}\int_I \int_M \left\{
g(\eta(t,x)) \left( u(t,\eta(t,x)), u(t,\eta(t,x))\right)\right. \nonumber\\
&\qquad\qquad + 
2\alpha^2 g(\eta(t,x))\left( \operatorname{Def}u(t,\eta(t,x)), 
\operatorname{Def} u(t,\eta(t,x))\right) \nonumber\\
&\qquad\qquad \left. + \lambda
g(\eta(t,x)) \left(
\nabla\left[ d(t,\eta(t,x)) \right],
\nabla\left[ d(t,\eta(t,x)) \right] \right) + 2 F (d)
\right\} \mu\, dt,
\label{averaged_action}
\end{align}
where we suppress the  explicit dependence of $u$ and $d$ on $\alpha$
and $\epsilon$.

Again, we see that (\ref{avg_lc.a}) arises as the first variation of the
action function $S^\alpha$ with respect to $\eta$, and the remaining equations
are identical to the original system (\ref{lc}).  Note, however, that now
$u$ is the {\it mean velocity}, and it is the {\it mean flow} $\eta$ which
is transporting the director field $d$.

\medskip

\noindent
{\bf Averaged Energy Law.}  For simplicity, we shall present the formulation
in the case that ${\partial M}=\emptyset$, although the more general case
follows in the same fashion as we presented above.
Following the notation of Section 
\ref{sec_variational}, we have the following basic averaged energy law:
\begin{align}
&{\frac{1}{2}}{\frac{d}{dt}} \left( |u|^2 + \alpha^2|\nabla u|^2
+ \lambda |\nabla d|^2 + 2\int_M F(d) \mu \right)\nonumber \\
&\qquad\qquad\qquad \le
-\nu\left( |\nabla u|^2 + \alpha^2 |\Delta u|^2\right) - \gamma \lambda
| \Delta d - f(d)|^2.
\label{avg_blaw}
\end{align}
From Lemma \ref{lemma3} and (\ref{avg_blaw}), it follows that there exists
$\bar t >0$ and some $\bar \rho_0 >0$ which are independent of the initial
data such that
$$ |u(t,\cdot)|_1^2 + |d|_1^2 \le \bar \rho_0 \ \ \ \ \forall t >\bar t.$$
We see that the averaged energy law is, in some sense, more natural than
the standard basic energy law  (\ref{energylaw}) since
the director field $d$ is no longer constrained to have one derivative 
greater regularity than the velocity of the fluid $u$:  both $u$ and $d$
now scale similarly.

Because of the a priori uniform bound of $u(t,\cdot)$ in $W^1$, it is very
easy to obtain an a priori bound for $u(t, \cdot)$ in $W^2$ and $d \in H^2(TM)$
when the $\operatorname{dim}(M)=3$.
We simply compute the sum
$$0=\langle (1-\alpha^2 \Delta) 
\text{(\ref{avg_lc.a})}, (1-\alpha^2) u\rangle + 
\langle \Delta \text{(\ref{avg_lc.c})}, \Delta d\rangle, $$
and use Lemma \ref{lemma3}.  Using similar estimates as above,  we obtain
an a priori energy estimate, in fact an absorbing set, in $W^2
\times H^2(TM)$, and by bootstrapping, we may easily obtain higher-order
a priori estimates.  Local well-posedness follows again from the contraction
mapping argument that we gave in Theorem \ref{thm1}, so we have the following

\begin{thm}\label{thm3}
For $n=2,3$, $s>1$, ${\partial M}= \emptyset$, $R=0$ and $u_0 \in W^s$, $d_0 \in
H^{s+1}(TM)$,
$$ u \in C^0([0,\infty], W^s), \ \ \ d\in C^0([0,\infty], H^{s+1}(TM)$$
are solutions to the system of equations (\ref{lc}).  Moreover, there
exists a compact global attractor for the system (\ref{lc})
in $W^{s-1} \times H^{s}(TM)$.
\end{thm}
It is not difficult to generalize this Theorem to manifolds with boundary
following the method in \cite{MS2}.

\section{Concluding remarks}
\noindent
{\bf Gradient flow versus damping.} 
We considered the $L^2$ gradient flow of the variation of the action function
$S(\eta,d)$ with
respect to $d$ in the director field equation (\ref{lc.c}).
In the liquid crystal literature, however, it is common to see a damped 
second-order equation for the director field (see \cite{DP} and
references therein), which in the context of our simplified system would
mean replacing $d_t$ with $\beta_1 d_{tt} + \beta_2 d_t$ for some constants
$\beta_1$ and $\beta_2$.  Of course, both types of equations have the
identical stationary solutions, but in terms of stability,
L. Simon's result \cite{Simon} guarantees that the damping term takes over.
As far as parabolic estimates are concerned, it is easy to treat either
type of equation, but we feel it is more natural to take the path of
steepest descent in relaxing the orientation towards its preferred 
configuration.

\medskip
\noindent
{\bf Lie advection versus parallel transport.}
This remark concerns the coupling term $\nabla_ud$
 in equation (\ref{lc.c}).  This term
arises by considering the time derivative of $(d \circ \eta)(t,x):=
d(t,\eta(t,x))$, where for
each $t$, $\eta(t,\cdot)$ is a volume-preserving diffeomorphism in
the topological group ${\mathcal D}_{\mu,D}^s$.  In 
group-theoretic language, this suggests that the action of 
${\mathcal D}_{\mu,D}^s$ on the vector space of director fields is on
the {\it right}.  The natural action of ${\mathcal D}_{\mu,D}^s$ on
the vector space of director fields, however,  is on the {\it left}, or by
push-forward:  instead of $d\circ \eta$, the natural action is
$\eta_* d := D\eta \cdot d \circ \eta^{-1}.$  Taking the time derivative
of $[\eta_* d](t,x)$  gives ${\mathcal L}_ud$, the Lie derivative of $d$
in the direction $u$.  The Lie derivative ${\mathcal L}_ud = \nabla_ud
- \nabla_du$, and this is the actual term which is present in the
Ericksen-Leslie model.  A nontrivial extension of our analysis is required to
analyze the system (\ref{lc}) with 
$\nabla_ud$ replaced by ${\mathcal L}_ud$, and this shall be the focus of
a future article.

\medskip \noindent
{\bf Other fluids models.}  Using our methodology, it is quite easy to study
a number of other fluids models.  For example, by replacing the Oseen-Frank
energy with the Landau-Lifshitz free energy $1/2\int_M A g(\nabla M, 
\nabla M) \mu$, where $M$ is the direction of magnetization in a cubic
ferromagnet, we can obtain an almost identical system of PDEs.  Similarly,
if we replace the vector $d$ in our action function $S$ with a scalar
field $\phi$, and replace the $L^2$ gradient flow in equation (\ref{lc.c})
with $H^{-1}$ gradient flow, we obtain a model of two-phase flow whose
interface moves via motion by mean curvature (see \cite{LS}).  
This model consists of a coupled Navier-Stokes Cahn-Hilliard system, 
where the interface is governed by surface tension.  

\medskip \noindent
{\bf The defect law in the limit as $\epsilon \rightarrow 0$.}  
We considered the GL penalization of the Oseen-Frank energy law so as
to obtain finite-energy minimizers, but we have yet to consider the limit
of our solutions as $\epsilon \rightarrow 0$.  It remains an open problem
to characterize the dynamical law of the GL vortices when coupled to
the Navier-Stokes motion.  Following the pioneering work in
\cite{BBH} and \cite{LX} on the dynamical law of the GL vortices, we expect that
the location of the jth vortex, $a_j$, will solve the distributional equation
$$ \partial_t a_j + \operatorname{div}(a_j u) = \frac{\delta W}{\delta a_j},$$
where $u$ simultaneously solves the Navier-Stokes equations, and
$W=\sum_{i\ne j} \log | x_i -x_j|$ is the renormalized energy.
We expect that a rigorous defect law will be much easier to obtain when
$u$ is instead a solution of the LANS equations, because in that case,
$u$ is uniformly in $H^1$ with respect to the penalization parameter
$\epsilon >0$.

\medskip
\begin{center}
{\sc Acknowledgments}
\end{center}
I am grateful to Chun Liu for introducing me to the beautiful theory of liquid
crystals.  I would also like to thank John Ball, Dick James, and Stefan 
M\"{u}ller 
for inviting me to the Oberwolfach conference Mathematical Continuum Mechanics 
where much of this work was initiated.
Research  was partially supported by the NSF-KDI grant ATM-98-73133
and the Alfred P. Sloan Foundation Research Fellowship.

\end{document}